\documentclass{article}

\usepackage{arxiv}

\usepackage{amsmath,amsfonts,bm}

\def\eqref#1{equation~\ref{#1}}

\def\1{\bm{1}}

\DeclareMathAlphabet{\mathsfit}{\encodingdefault}{\sfdefault}{m}{sl}
\SetMathAlphabet{\mathsfit}{bold}{\encodingdefault}{\sfdefault}{bx}{n}

\newcommand{\R}{\mathbb{R}}

\usepackage{amsmath,amsfonts,bm}

\def\eqref#1{equation~\ref{#1}}

\def\1{\bm{1}}

\DeclareMathAlphabet{\mathsfit}{\encodingdefault}{\sfdefault}{m}{sl}
\SetMathAlphabet{\mathsfit}{bold}{\encodingdefault}{\sfdefault}{bx}{n}

\usepackage{algorithm}
\usepackage{algorithmic}
\usepackage{amsthm}
\usepackage{xcolor}
\usepackage{subcaption}

\usepackage[utf8]{inputenc} 
\usepackage[T1]{fontenc}    
\usepackage{hyperref}       
\usepackage{url}            
\usepackage{booktabs}       
\usepackage{amsfonts}       
\usepackage{nicefrac}       
\usepackage{microtype}      
\usepackage{lipsum}
\usepackage{graphicx}
\usepackage{cleveref}
\usepackage{verbatim}
\usepackage{float}
\usepackage{multicol}
\usepackage{soul}
\usepackage{fancyhdr}

\title{Sparse Random-Feature Neural Networks with Krylov-Based SVD for Singularly Perturbed ODE}

\author{
 Kevin Kurian Thomas Vaidyan$,$ Siddharth Rout\\
  Institute of Applied Mathematics\\
  University of British Columbia\\
  Vancouver, BC V6T 1Z4 \\
}

\begin{document}
\maketitle
\begin{abstract}
Random-feature neural networks (RFNNs), including architectures with fixed hidden layers and analytically determined output weights, offer fast training but often suffer from issues due to dense representations of the hidden layer activation. Their reliance on dense feature mappings and least squares solvers can limit scalability and numerical stability, particularly for high-dimensional or stiff systems. Specifically, the activation matrix is observed to be low-rank and extremely ill-conditioned. In this work, we propose a sparse framework for RFNNs that integrates structured sparsity into the hidden layer activations that increases the rank and employs Sparse Singular Value Decomposition (sSVD) for solving the resulting linear least squares problem scalably and efficiently while catering to the bad condition number. We explore the theory behind Lanczos-Golub-Kahan Bidiagonalization technique for sparse SVD and conduct some experiments to identify some limitations and justify the requirement for orthogonalization step in our application. Then, we demonstrate that the proposed method maintains or improves solution accuracy for solving the benchmark \textbf{one-dimensional steady convection-diffusion equations case having stronger advection}, while achieving substantial gains in training efficiency and robustness compared to standard dense implementations.

\end{abstract}

\section{Introduction}

Neural networks have demonstrated remarkable versatility and effectiveness across a broad range of applications, including images and audios \cite{8632885,9451544,peebles2023scalable}, natural language processing \cite{peebles2023scalable,chang2024survey,thirunavukarasu2023large}, and complex scientific modelling \cite{thiyagalingam2022scientific,cuomo2022scientific, pathak2022fourcastnet}. However, neural networks have often failed to cater the scientific community where precision is of utmost importance. For an instance, it is a challenge for neural networks to learn some simple low dimensional mathematical or physical symbolic functions with a misfit in the order of machine precision \cite{liu2024kan}. Problems arise when the data fitting is not accurate enough. Such problems often appear when simulations are strongly sensitive to the look-up variables like in case of chaotic dynamics in turbulent fluid dynamics in weather models, combustion, etc. When it comes to accurate prediction for complex scientific applications, scalability, accuracy and simplicity of the machine learning approximator are of utmost importance. 

Random-feature neural networks (RFNNs \cite{huang2006extreme, misiakiewicz2023six}), also called extreme learning machines (ELMs), including architectures with fixed hidden layers and analytically determined output weights, have emerged as efficient alternatives to conventional gradient-based learning methods. By decoupling feature generation from parameter optimization, these models enable rapid training through closed-form solutions, typically framed as least squares problems. Such characteristics make RFNNs particularly attractive for applications requiring fast inference and reduced training complexity. This construction transforms nonlinear learning into a convex optimization problem over the output weights, avoiding the non-convexity. Amazingly, despite such simplifications, the theory of universal approximation remains valid\cite{Huang2006UniversalAU, HUANG20073056}. The beauty of this method, which makes it special, is that the loss function could be represented as a set of equations that are nonlinear in terms of input variables and, however, linear in terms of trainable weights. Random-feature models are closely related to kernel methods, as they approximate kernel-induced feature maps while enabling efficient large-scale computation \cite[see section 9.5, page number 271]{bach2024learning}. There have been many models in this class of approximators. Some popular and impactful models include echo-state network \cite{jaeger2004harnessing}, extreme learning machine \cite{huang2006extreme} and randomized kernel machine \cite{rahimi2007random}.

Despite their simplicity and favourable optimization properties, the performance of random-feature neural networks depends critically on the quality of the induced feature representations \cite{rudi2017generalization, mei2022generalization}. However, despite these advantages, their reliance on dense feature representations often leads to scalability challenges, and increased memory consumption, especially in high-dimensional settings. Again, ill-conditioning, and rank deficiency pose problems from a numerical optimization standpoint. These challenges motivate the development of robust feature processing and efficient numerical methods to improve stability and scalability without sacrificing predictive accuracy. 

Particularly, in the domain of numerical partial differential equations, physics-Informed Extreme Learning Machines (PIELMs \cite{dwivedi2020physics}), a specialized class of RNNs, extend this paradigm to scientific machine learning by embedding physical laws directly into the learning process. Similar to physics-informed neural networks (PINNs \cite{raissi2019physics, karniadakis2021physics, cuomo2022scientific}), PIELMs enforce governing equations—such as partial differential equations (PDEs)—and boundary conditions within the training objective. We considered that PINNs are, essentially, complex stiff optimization problems as the algorithm requires expensive optimizers like L-BFGS \cite{raissi2019physics, wang2021understanding, mckay2025near, cao2025analysis} or techniques like domain-decomposition \cite{rout2019numerical, jagtap2020extended}, preconditioning \cite{kopanivcakova2024enhancing, liu2024preconditioning, song2025matrix, ahmad2025preconditioned}, multi-staging \cite{wang2024multi, yuan2024physics, khadijeh2025multistage}. This seriously limits the use of large neural networks. Unlike PINNs, however, PIELMs avoid iterative gradient-based optimization by leveraging randomized hidden representations and computing output weights analytically. This results in significantly reduced training times while maintaining the ability to solve complex physical systems and while allowing easy analysis of the learning algorithm. Consequently, PIELMs have gained attention as a computationally efficient framework for solving linear forward and inverse problems in physics-based modeling like in heat transfer \cite{ren2025physics}, linear elasticity \cite{wang2025physics},  and Eigenmode analysis \cite{mishra2025eig, huang2025physics}.

Despite their promises, PIELMs inherit key limitations from their underlying randomized architectures. In particular, the dense construction of hidden layer activations and the associated least squares systems can become computationally expensive and numerically ill-conditioned as the number of collocation points or hidden neurons increases\cite{rout2026fast}. These issues are exacerbated in large-scale or high-resolution simulations, where both memory usage and computational cost grow rapidly.

To address these challenges, recent efforts have explored the incorporation of sparsity and advanced linear algebra techniques within random-feature neural network frameworks \cite{liu2011selm, luo2013sparse, bai2014sparse, yu2017sparse}. Sparse representations have long been recognized for their ability to reduce computational complexity, improve memory efficiency, and enhance numerical stability in large-scale systems. In parallel, sparse matrix factorization methods—such as sparse singular value decomposition (SVD)—offer robust tools for solving ill-conditioned least squares problems while exploiting underlying structural sparsity. Motivated by these developments, this work investigates the integration of sparsity into the PIELM framework. By introducing sparse feature representations and leveraging sparse SVD for the resulting least squares formulation, we aim to improve the scalability, efficiency, and robustness of PIELMs for solving complex physical systems. The proposed approach seeks to retain the computational advantages of random-feature neural networks while addressing their key limitations in dense settings.

Here are our contributions explained in this work:
\begin{itemize}
\item A novel algorithm for a robust random-feature neural network with a linear readout layer by combining Gaussian filtration and sparse matrix computation.

\item A technique using Gaussian filtration to improve the effective rank of random feature representations, enhancing learnability.

\item Exploited sparsity-aware SVD algorithm that significantly reduce memory usage and computational complexity for training.

\item Studied the theory of Golub-Kahan Bidiagonalization and the need for full-orthogonalization using Gram-Schmidt process. 

\item Demonstrated through experiments improved robustness and substantial computational gains compared to dense random-feature baselines.
\end{itemize}

\section{Random-Feature Neural Networks}
A random-feature neural network $y$ with a linear readout layer whose trainable weights are \(\beta\), for activation function $\phi$ and random matrices $W$ and $b$ can be written mathematically as:\\
\begin{equation}
\label{eqn:elmarchitecture}
y(\mathbf{x}) = \phi(\mathbf{W} \mathbf{x} + \mathbf{b}) \mathbf{\beta}.
\end{equation}


\textbf{Key Property:} Nonlinear in $\mathbf{x}$ while linear in trainable weights $\mathbf{\beta}$.

\textbf{Training via least squares:}
\[
\min_{\mathbf{\beta}} \; \| \Phi \mathbf{\beta} - Y \|^2,
\]

where $\Phi$ is the activation matrix corresponding to the input matrix $X = \{x_1, ..., x_N\}$ and output matrix $Y = \{y_1, ..., y_N\}$. \\

\textbf{Closed-Form Solution:}
\[
\mathbf{\beta}^* = (\Phi^T \Phi)^{-1} \Phi^T Y
\quad \text{(or pseudoinverse)}.
\]

\subsection{Solving Differential Equations using Neural Networks}
Let us say for any differential equation $\mathcal{N}$ defined by \begin{equation}
\label{eq:1}
\mathcal{N}[u](\mathbf{x}) = 0,
\end{equation}
with constraints at $\mathbf{X}_u$ we have 
\begin{equation}
\label{eq:2}
u(\mathbf{X}_u) = \mathbf{u}_o.\\
\end{equation} 

The algorithms for solving the PDE defined by \cref{eq:1} and \cref{eq:2} using PINN and PIELM are mentioned in the subsequent subsections.

\subsection{Physics-Informed Neural Networks}
The \cref{pinn_algo} shows the pseudocode for PINNs as per the vanilla settings in \cite{raissi2019physics}. The method uses training data points $(\mathbf{X}_u, \mathbf{u})$ to fit the neural network with trainable parameters $\theta$ denoted by $u_\theta(\mathbf{x})$, and minimize the residual of $\mathcal{N}[u](\mathbf{x})$ at collocation points $\mathbf{X}_f$.

\begin{algorithm}[!h]
\caption{PINN}
\begin{algorithmic}[1]
\label{pinn_algo}

\STATE \textbf{Initialization}
Initialize neural network parameters $\theta$ (weights and biases)

\STATE \textbf{Define the Neural Network}
Construct the neural network $u_\theta(\mathbf{x})$ with parameters $\theta$

\STATE \textbf{Formulate the Loss Functions}
Define data fitting loss:
\begin{equation}
\label{eq: datafit}
\mathcal{L}_{\text{data}} = \frac{1}{N_u} \sum_{i=1}^{N_u} \left| u_\theta(\mathbf{x}_u^i) - u^i \right|^2.
\end{equation}

Define physics-informed loss using the governing differential equation $\mathcal{N}[u] = 0$:
\begin{equation}
\label{eq: pdefit}
\mathcal{L}_{\text{phys}} = \frac{1}{N_f} \sum_{j=1}^{N_f} \left| \mathcal{N}[u_\theta](\mathbf{x}_f^j) \right|^2.
\end{equation}
Combine the losses:
\begin{equation}
\label{eq: multiobj}
\mathcal{L} = \mathcal{L}_{\text{data}} + \mathcal{L}_{\text{phys}}.
\end{equation}

\STATE \textbf{Training}
Minimize the total loss $\mathcal{L}$ with respect to the network parameters $\theta$:
\begin{equation}
\theta^* = \arg\min_\theta \mathcal{L}(\theta).
\end{equation}

\STATE \textbf{Inference}
Use the trained network $u_{\theta^*}$ to make predictions for new inputs $\mathbf{X}_{\text{new}}$:
\begin{equation}
\mathbf{\hat{u}} = u_{\theta^*}(\mathbf{X}_{\text{new}}).
\end{equation}

\end{algorithmic}
\end{algorithm}

\subsection{Physics-Informed Extreme Learning Machine (PI-ELM) Algorithm}
The \cref{algo_pielm} shows the pseudocode for PI-ELM. Similar to the problem definition in PINNs, this method also uses training data points $(\mathbf{X}_u, \mathbf{u})$ to fit the neural network with trainable parameters $\beta$ denoted by $u_\beta(\mathbf{x})$ and minimise the residual of $\mathcal{N}[u](\mathbf{x})$ at collocation points $\mathbf{X}_f$. Trained output weights $\mathbf{\beta}$ are obtained by solving the system of linear equations.

\begin{algorithm}[!h]
\caption{PI-ELM}
\begin{algorithmic}[1]
\label{algo_pielm}

\STATE \textbf{Initialization}
Randomly initialize input weights $\mathbf{W} \in \mathbb{R}^{L \times n}$ and biases $\mathbf{b} \in \mathbb{R}^L$.
Define activation function $\phi(\cdot)$

\STATE \textbf{Formulate Composite Loss Function}
Define data fitting loss: 
\begin{equation}
\mathcal{L}_{\text{data}} = \sum_{i=1}^{N_u} \left| u(\mathbf{x}_u^i) - u^i \right|^2.
\end{equation}
Define physics-informed loss:
\begin{equation}
\mathcal{L}_{\text{phys}} = \sum_{i=1}^{N_f} \left| \mathcal{N}[u(\mathbf{x}_f^i)] \right|^2.
\end{equation}
Combine the losses:
\begin{equation}
\mathcal{L} = \mathcal{L}_{\text{data}} + \mathcal{L}_{\text{phys}}.
\end{equation}

\STATE \textbf{Generate the System of Equations}
The equation for the neural network $u_\beta(\mathbf{x})$ is
\begin{equation}
\label{eq:elm}
u_\beta(\mathbf{x}) = \phi(\mathbf{W} \mathbf{x} + \mathbf{b}) \mathbf{\beta}.
\end{equation}

The effective system of linear equations in $\mathcal{L}$ could be rewritten in the form:
\begin{equation}
\label{eq:PIELM_Inversion}
\mathbf{H} \mathbf{\beta}  = \mathbf{T}.
\end{equation}

\STATE \textbf{Training}
Solve for the output weights using the Moore-Penrose pseudoinverse of $\mathbf{H}$ denoted by $\mathbf{H}^+$ which is essentially the least square solution:
\begin{equation}
\mathbf{\beta} = \mathbf{H}^+ \mathbf{T}.
\end{equation}

\STATE \textbf{Inference}
Use the trained PI-ELM model to predict outputs for new inputs $\mathbf{X}_{\text{new}}$:
\begin{equation}
\mathbf{\hat{Y}} = \phi(\mathbf{W} \mathbf{X}_{\text{new}} + \mathbf{b}) \mathbf{\beta}.
\end{equation}

\end{algorithmic}
\end{algorithm}

\section{Sparse Singular Value Decomposition}
Singular Value Decomposition (SVD) is one of the most fundamental tools in numerical linear algebra, underpinning a wide range of applications in scientific computing, machine learning, signal processing, and data analysis. Given a matrix 
\( A \in \mathbb{R}^{m \times n},\) 
the SVD provides a factorization $A = U \Sigma V^\top$. 

However, classical SVD algorithms are computationally expensive. This is especially true when the matrix $A$ is large and sparse, which is often the case in real-world applications such as natural language processing and recommender systems.

In these settings, applying standard dense SVD methods is both computationally inefficient and memory-intensive, as they fail to exploit the underlying sparsity structure. This motivates the development of \emph{sparse SVD methods}, which aim to compute a subset of singular values and vectors—typically the largest or smallest—while preserving sparsity and avoiding unnecessary operations on zero entries.

In this report, we explore the Lanczos bidiagonalization algorithm, a powerful iterative method for computing the SVD of large sparse matrices and its restarted version. We will also evaluate the performance of these algorithms on various sparse matrices and test it on a neural PDE problem.

\subsection{A motivating example: Eigenvalue problems}
We begin with a brief review of the Lanczos method for computing eigenvalues and eigenvectors of large sparse symmetric matrices \cite{golub2013matrix}. The Lanczos method is an iterative algorithm that constructs a sequence of orthogonal vectors and tridiagonal matrices, which can be used to approximate the eigenvalues and eigenvectors of the original matrix. Consider a symmetric matrix $A \in \mathbb{R}^{n \times n}$ and an initial vector $v_1$ with $\|v_1\| = 1$. The Lanczos method generates a sequence of vectors $\{v_k\}$ and scalars $\{\alpha_k, \beta_k\}$ such that
\[Av_k = V_k T_k + \beta_{k+1} v_{k+1} e_k^\top,\]
where $V_k = [v_1, v_2, \ldots, v_k]$ is an orthonormal basis for the Krylov subspace 
\[ K_k(A; v_1) = \text{span}\{v_1, Av_1, A^2v_1, \ldots, A^{k-1}v_1\} \]
and $T_k$ is a tridiagonal matrix with diagonal entries $\alpha_k$ and off-diagonal entries $\beta_k$. In class, we discussed how the eigenvalues of $T_k$ can be used to approximate the eigenvalues of $A$, and how the corresponding eigenvectors can be approximated using the Lanczos vectors. Suppose $(\theta, y)$ is an eigenpair of $T_k$ and define $u = V_k y$. Now notice,
\begin{align*}
    Au - \theta u &= AV_k y - \theta V_k y \\
    &= V_k T_k y + \beta_{k+1} v_{k+1} e_k^\top y - \theta V_k y \\
    &= V_k ( T_k y - \theta y) + \beta_{k+1} v_{k+1} e_k^\top y \\
    &= \beta_{k+1} v_{k+1} e_k^\top y.
\end{align*}
Thus, the residual norm $\|Au - \theta u\|$ is equal to $\beta_{k+1} |e_k^\top y|$. If $\beta_{k +1}$ is small and/or $|e_k^\top y|$ is small, then $u$ is a good approximation to an eigenvector of $A$ corresponding to the eigenvalue $\theta$. We would also like to point out that if $\beta_{k+1} = 0$, (when the Krylov subspace is invariant and not expanding) then $u$ is an exact eigenvector of $A$. 

It is known that the Lanczos method can converge very quickly to the largest and smallest eigenvalues of $A$. As a summary, we compute the eigenvectors-$y_k$ and eigenvalues-$\theta_k$ of $T_k$ and the approximate eigenvalues of $A$ are given by $\theta_k$. The corresponding eigenvectors of $A$ are given by $u_k = V_k y_k$. 

\subsection{A naive approach to sparse SVD}
One naive approach that can be used to compute the SVD of a large sparse matrix $A \in \R^{m \times n}$ is to notice that
\[\sigma_i(A) = \sigma_i(A^\top A)\]
So one algorithm is
\begin{algorithm}
\caption{Naive SVD algorithm}
\begin{algorithmic}[1]
\STATE Compute $B = A^\top A$.
\STATE Compute the eigenvalues $\lambda_i$ and eigenvectors $v_i$ of $B$ using the Lanczos method.
\STATE Compute the singular values $\sigma_i = \sqrt{\lambda_i}$, the right singular vectors $v_i$ and the left singular vectors $u_i = Av_i / \sigma_i$.
\end{algorithmic}
\end{algorithm}
This algorithm is not efficient for large sparse matrices because:
\begin{itemize}
    \item The matrix $B = A^\top A$ is typically dense, even if $A$ is sparse. This means that computing $B$ and storing it in memory can be very expensive.
    \item The eigenvalues of $B$ are the squares of the singular values of $A$. If we are only interested in a few of the largest singular values, we may end up computing many eigenvalues of $B$ that are not relevant to our problem.
    \item The condition number of $B$ is the square of the condition number of $A$. This can lead to numerical instability.
\end{itemize}

Hence, we need more sophisticated algorithms that can compute the SVD of large sparse matrices without explicitly forming $A^\top A$. In the next section, we will discuss the Golub-Kahan bidiagonalization algorithm, which is an efficient method for computing the SVD of large sparse matrices. 
We will also explain the connection between the Lanczos method and the Golub-Kahan bidiagonalization algorithm.

\subsection{Golub-Kahan bidiagonalization algorithm}
The Golub-Kahan bidiagonalization algorithm is an iterative method for decomposing a matrix $A \in \R^{m \times n}$ into a product of two orthonormal matrices ($U_k \in \R^{m \times k}, V_k \in \R^{n \times k}$) and a bidiagonal matrix $B \in \R^{k \times k}$. 
Given an initial unit vector $v_1 \in \mathbb{R}^n$, the Golub--Kahan process constructs two orthonormal sequences

\[
\{u_1, \dots, u_k\} \subset \mathbb{R}^m, \quad \{v_1, \dots, v_k\} \subset \mathbb{R}^n,
\]
along with scalars $\{\alpha_j\}, \{\beta_j\}$ satisfying the recurrences
\begin{align*}
    &A v_j = \beta_{j-1} u_{j-1} + \alpha_j u_j, \\
    &A^T u_j = \alpha_j v_j + \beta_j v_{j+1},
\end{align*}
for $j = 1, \dots, k$, with initialization $\beta_0 = 0$, $\alpha_1 = \|A v_1\|$ and $u_1 = A v_1 / \alpha_1$. 

By construction, the vectors $\{u_j\}$ and $\{v_j\}$ are orthonormal, and the scalars $\alpha_j, \beta_j \geq 0$. These vectors form a orthonormal basis for the Krylov subspaces:

\begin{align}
    &\text{span}\{v_1, \dots, v_k\} = \mathcal{K}_k(A^TA, v_1) = \text{span}\{v_1, (A^TA)v_1, \ldots, (A^TA)^{k-1} v_1\},\\
    &\text{span}\{u_1, \dots, u_k\} = \mathcal{K}_k(AA^T, u_1) = \text{span}\{u_1, (AA^T)u_1, \ldots, (AA^T)^{k-1} u_1\}.
\end{align}

\subsubsection{Matrix Formulation}

Let
\[
U_k = [u_1, \dots, u_k] \in \mathbb{R}^{m \times k}, \quad
V_k = [v_1, \dots, v_k] \in \mathbb{R}^{n \times k}.
\]
Then the recurrences can be written compactly as
\begin{align}
    A V_k &= U_k B_k, \\
    A^\top U_k &= V_k B_k^\top + \beta_{k} v_{k+1} e_k^\top,
\end{align}
where $B_k \in \mathbb{R}^{k \times k}$ is an upper bidiagonal matrix of the form
\[
B_k =
\begin{bmatrix}
\alpha_1 & \beta_2 & 0 & \cdots & 0 \\
0 & \alpha_2 & \beta_3 & \ddots & \vdots \\
0 & 0 & \alpha_3 & \ddots & 0 \\
\vdots & \ddots & \ddots & \ddots & \beta_k \\
0 & \cdots & 0 & 0 & \alpha_k
\end{bmatrix}.
\]
If we multiply the first equation by $U_k^\top$, (or equivalently, the second equation by $V_k^\top$) we get
\begin{equation}
    U_k^\top A V_k = B_k
\end{equation}
The matrix $B_k$ is a low-dimensional representation of $A$ in the Krylov subspaces generated by $A^\top A$ and $A A^\top$. 

\subsubsection{Computing singular values and vectors}

The singular values of $B_k$ serve as approximations (called \emph{Ritz singular values}) to the singular values of $A$. If
\[
B_k = \widetilde{U} \Sigma \widetilde{V}^\top
\]
is the SVD of $B_k$, then notice that, 
\begin{align*}
    A (V_k \widetilde{V}) &= U_k B_k \widetilde{V} \\
    &= U_k \widetilde{U} \Sigma.
\end{align*}

Similarly, we can also show that
\begin{align*}
    A^\top (U_k \widetilde{U}) &= V_k B_k^\top \widetilde{U} + \beta_{k} v_{k+1} e_k^\top \widetilde{U} \\
    &= V_k \widetilde{V} \Sigma + \beta_{k} v_{k+1} e_k^\top \widetilde{U}.
\end{align*}

Thus, the Ritz singular values $\Sigma$ are approximations to the singular values of $A$, and the Ritz singular vectors $U_k \widetilde{U}$ and $V_k \widetilde{V}$ are approximations to the left and right singular vectors of $A$, respectively.
Also, notice that the residual norm $\|A^\top (U_k \widetilde{U}) - (V_k \widetilde{V}) \Sigma\|$ is equal to $\beta_{k} \|e_k^\top \widetilde{U}\|$. If $\beta_{k}$ is small and/or $\|e_k^\top \widetilde{U}\|$ is small, then the Ritz singular vectors are good approximations to the singular vectors of $A$ corresponding to the Ritz singular values.

\subsection{Relation to Lanczos diagonalization}
In this section, we will explain the connection between the Lanczos method and the Golub-Kahan bidiagonalization algorithm that we had failed to address during our midterm project. 

We claim that the Lanczos method applied to $A^\top A$ produces the same sequence of vectors $\{v_j\}$ and scalars $\{\alpha_j, \beta_j\}$ as the Golub-Kahan bidiagonalization algorithm applied to $A$. Starting from the first recurrence of the Golub-Kahan algorithm, we have
\begin{align*}
    A V_k &= U_k B_k\\
\end{align*}
Multiplying both sides by $A^\top$, we get
\begin{align*}
    A^\top A V_k &= A^\top U_k B_k \\
    &= V_k B_k^\top B_k + \beta_{k} v_{k+1} e_k^\top \alpha_k \tag{Last entry is $\alpha_k$}\\
    &= V_k T_k + \alpha_k \beta_{k} v_{k+1} e_k^\top \\
\end{align*}
where $T_k = B_k^\top B_k$ is a tridiagonal matrix. This is exactly the Lanczos recurrence for $A^\top A$ with the same vectors $\{v_j\}$ and scalars $\{\alpha_j, \beta_j\}$. \\

We provide pseudocode for the Golub-Kahan bidiagonalization algorithm below:
\begin{algorithm}
\caption{Golub--Kahan Bidiagonalization}
\begin{algorithmic}[1]
\STATE Choose $v_1 \in \mathbb{R}^n$ with $\|v_1\|_2 = 1$
\STATE Set $\beta_0 = 0$ and $u_0 = 0$
\FOR{$j = 1,2,\dots,k$}
    \STATE $u_j = A v_j - \beta_{j-1} u_{j-1}$
    \STATE $\alpha_j = \|u_j\|_2$
    \STATE $u_j = u_j / \alpha_j$
    \STATE $v_{j+1} = A^\top u_j - \alpha_j v_j$
    \STATE $\beta_j = \|v_{j+1}\|_2$
    \STATE $v_{j+1} = v_{j+1} / \beta_j$
\ENDFOR
\end{algorithmic}
\end{algorithm}

\subsection{Dealing with loss of orthogonality}
In the idealized version of the Golub-Kahan algorithm, the vectors $\{u_j\}$ and $\{v_j\}$ are orthonormal. That is
\[V_k^\top V_k = I_k \quad \text{and} \quad U_k^\top U_k = I_k.\]
However, in practice, due to finite precision arithmetic, the vectors can lose orthogonality. This can lead to numerical instability and inaccurate approximations of the singular values and vectors. 
There are several different approaches to mitigate the loss of orthogonality. We discuss a few below:

\subsubsection{Full orthogonalization}
The bidiagonalization procedure produces the set of left singular vectors $\{u_j\}$ and right singular vectors $\{v_j\}$. To maintain orthogonality, we can perform a full reorthogonalization step after each iteration. This involves orthogonalizing the new vector against all previously computed vectors in the respective set. 
As a reminder, this is very similar to the initial CG algorithm we developed in class where we performed a full A-conjugate gram-schmidt step to find A-conjugate directions. The algorithm is provided below:
\begin{algorithm}[H]
\caption{Lanczos Bidiagonalization with Full Orthogonalization}
\begin{algorithmic}[1]
\STATE Choose a unit vector $v_1 \in \mathbb{R}^n$, $\|v_1\| = 1$
\FOR{$j = 1,2,\dots,k$}

    \STATE $u_j = A v_j$
    \FOR{$i = 1,2,\dots,j-1$}
        \STATE $\gamma = u_i^\top u_j$
        \STATE $u_j = u_j - \gamma u_i$
    \ENDFOR
    \STATE $\alpha_j = \|u_j\|$
    \STATE $u_j = u_j / \alpha_j$

    \STATE $v_{j+1} = A^\top u_j$
    \FOR{$i = 1,2,\dots,j$}
        \STATE $\gamma = v_i^\top v_{j+1}$
        \STATE $v_{j+1} = v_{j+1} - \gamma v_i$
    \ENDFOR
    \STATE $\beta_j = \|v_{j+1}\|$
    \STATE $v_{j+1} = v_{j+1} / \beta_j$

\ENDFOR
\end{algorithmic}
\end{algorithm}

Lines 4-8 and 11-14 perform the full gram-schmidt orthogonalization step to maintain the orthogonality of the vectors. However, this method has the same drawbacks that we discussed 
in class for the original CG algorithm we developed. The computational cost of the full orthogonalization step is $O(j)$ for the $j$-th iteration, which can become expensive as $j$ increases.

\subsubsection{One sided orthogonalization}
Another approach is to perform orthogonalization on only one set of vectors. This technique is proposed by Simon and Zha \cite{simon2000low} and is called one-sided orthogonalization. The main idea behind this method is that the level of orthogonality of the left and right lanczos vectors are pretty similar. If $V_k^*, U_k^*$ are the finite precision versions of $V_k, U_k$ respectively, then we can show that
$\|I - V_k^* V_k\|$ and $\|I - U_k^* U_k\|$ are of the same order except when $B_k$ becomes ill conditioned. In their method, they perform full orthogonalization on the right singular vectors $\{v_j\}$ and the left singular vectors $\{u_j\}$ are generated like in the original Golub-Kahan algorithm. The algorithm is provided below:

\begin{algorithm}
\caption{One-Sided Lanczos Bidiagonalization}
\begin{algorithmic}[1]
\STATE Choose a unit vector $v_1 \in \mathbb{R}^n$, $\|v_1\| = 1$
\STATE Set $\beta_0 = 0$ and $u_0 = 0$

\FOR{$j = 1,2,\dots,k$}

    \STATE $u_j = A v_j - \beta_{j-1} u_{j-1}$
    \STATE $\alpha_j = \|u_j\|$
    \STATE $u_j = u_j / \alpha_j$

    \STATE $v_{j+1} = A^\top u_j$
    \FOR{$i = 1,2,\dots,j$}
        \STATE $\gamma = v_i^\top v_{j+1}$
        \STATE $v_{j+1} = v_{j+1} - \gamma v_i$
    \ENDFOR
    \STATE $\beta_j = \|v_{j+1}\|$
    \STATE $v_{j+1} = v_{j+1} / \beta_j$

\ENDFOR
\end{algorithmic}
\end{algorithm}

The main drawback of this method is that it is sensitive to the condition number of $B_k$. If $B_k$ becomes ill-conditioned, then we would have to revert to the 2 sided orthogonalization method to maintain the accuracy.


\section{Proposed Method}
\subsection{Functional architecture}
The architecture we propose is a layer of encoding after the layer of linear combination of weights to the inputs. Following the procedure of algorithm \cref{algo_pielm}, we make a change to the equation \cref{eq:elm} following the same set of definitions and notations. The replacement equation is stated as:
\begin{equation}
\label{encoded_ELM}
u_\beta(\mathbf{x}) = \phi(\mathbf{W} \mathbf{x} \odot E(\mathbf{x}, \mu) + \mathbf{b}) \mathbf{\beta},
\end{equation}
where $E$ is the encoding function, $\odot$ is a Hadamard product and $\mu$ is a hyperparameter.

\subsubsection{Shifted Gaussian Encoding} 
In this work, the encoding function $E$ is defined as:
\begin{equation}
E(x) = e^{-\frac{(x - \mu)^2}{d}},    
\end{equation}
where $E$ is a Gaussian Kernel or Radial Basis Function (RBF), d is called the filter width, $\mu$ is a set {$\mu_i = \frac{i}{L} | i \in [0, 1,..,L]$}, and $L$ is one less than the number of hidden nodes. 

\subsubsection{Sparse SVD-Based Pseudoinverse}

Let $A \in \mathbb{R}^{m \times n}$ admit the singular value decomposition (SVD)
\begin{equation}
A = U \Sigma V^\top,    
\end{equation}
where $U \in \mathbb{R}^{m \times m}$ and $V \in \mathbb{R}^{n \times n}$ are orthogonal matrices, and $\Sigma \in \mathbb{R}^{m \times n}$ is a diagonal matrix containing the singular values $\{\sigma_i\}_{i=1}^r$ with $r = \operatorname{rank}(A)$.

The Moore--Penrose pseudoinverse of $A$ is given by
\begin{equation}
A^{\dagger} = V \Sigma^{\dagger} U^\top,
\end{equation}
where $\Sigma^{\dagger}$ is obtained by taking the reciprocal of the nonzero singular values:
\begin{equation}
\Sigma^{\dagger} = \operatorname{diag}\left(\frac{1}{\sigma_1}, \frac{1}{\sigma_2}, \dots, \frac{1}{\sigma_r}, 0, \dots \right).
\end{equation}

Using this formulation, the solution to the least squares problem
\begin{equation}
\min_{\beta} \|A \beta - Y\|_2^2
\end{equation}
is given by
\begin{equation}
\beta = A^{\dagger} Y = V \Sigma^{\dagger} U^\top Y.
\end{equation}

\paragraph{Truncated SVD.}
In practice, small singular values may lead to numerical instability. A truncated SVD can be used:
\[
\sigma_i^{-1} =
\begin{cases}
\frac{1}{\sigma_i}, & \sigma_i > \varepsilon, \\
0, & \text{otherwise},
\end{cases}
\]
which acts as a regularization mechanism and improves stability for ill-conditioned problems.

\paragraph{Sparse Setting.}
For large sparse matrices, forming the full SVD is computationally expensive. Instead, Krylov subspace methods such as Golub--Kahan bidiagonalization are employed to compute a low-rank approximation:
\begin{equation}
A \approx U_k B_k V_k^\top,
\end{equation}
leading to the approximate pseudoinverse
\begin{equation} \label{equ:sparsePseudoINV}
A^{\dagger} \approx V_k B_k^{-1} U_k^\top.
\end{equation}

\section{Experiments}

\subsection{Understanding the problem: Solving Non-smooth ODE (Steady 1D Convection-Diffusion)} \label{problem}

\textbf{Governing Equation:}
\begin{equation}
    u \frac{d\phi}{dx} = D \frac{d^2 \phi}{dx^2}    
\end{equation} 

\textbf{Boundary Conditions:} 
\begin{equation}
   \phi(0) = \phi_0, \quad \phi(L) = \phi_L 
\end{equation}

\textbf{Exact Solution:} 
\begin{equation}
    \phi(x) = \phi_0 + (\phi_L - \phi_0) \frac{e^{\frac{u}{D}x} - 1}{e^{\frac{u}{D}L} - 1} 
\end{equation}\\

In this work, we take $L = 1$, $\phi_0 = 0$,  $\phi_L = 1$ and $Pe = \frac{u}{D} = -1e3$.\\

\textbf{Key Behavior:} \begin{itemize} 
\item $Pe \gg 1$: Convection-dominated $\rightarrow$ Ultra Sharp Boundary layer 
\end{itemize} 

\begin{figure}[!t]
    \centering
\begin{minipage}{0.49\textwidth}
\includegraphics[width=\textwidth]{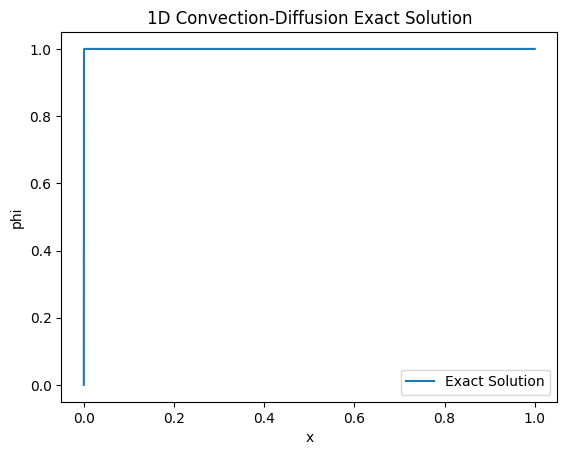}
\end{minipage}
\hfill
\begin{minipage}{0.49\textwidth}
\includegraphics[width=\textwidth]{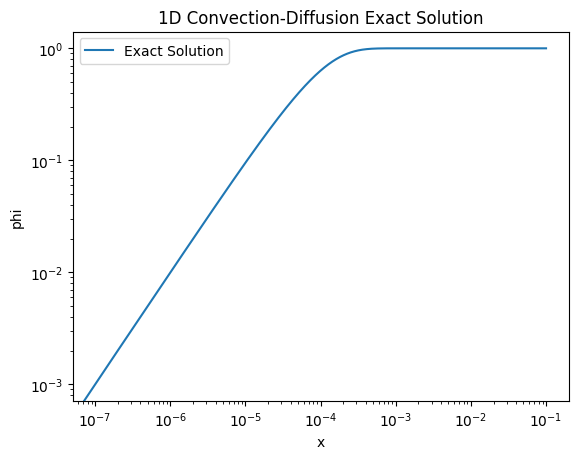}
\end{minipage}
    \caption{True solution for $Pe = -1e4$.}
    \label{fig:placeholder}
\end{figure}

\newpage
\subsubsection{Comparison of activation matrices (for reference)}
\begin{figure}[!b]
\centering
\begin{minipage}{0.25\textwidth}
\includegraphics[width=\textwidth]{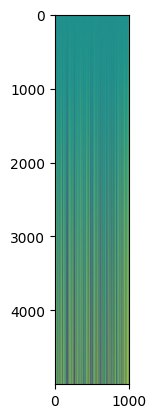}
PIELM
\end{minipage}
\begin{minipage}{0.25\textwidth}
\includegraphics[width=\textwidth]{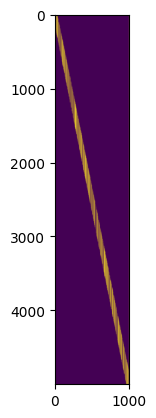}
Sparse-PIELM (ours)
\end{minipage}
\caption{Visualization of activation matrices.}
\label{fig:matrices}
\end{figure}

\Cref{tab:matrices} shows the properties of the matrix formed in the linear least squares setup. \Cref{fig:matrices} shows the illustration of matrix items coloured by the magnitudes. Our formulation makes the matrix sparse, while the sparsity is controlled by the width sof Gaussian Kernel. Decreasing the kernel width makes the matrix very sparse and increases the rank.

\subsubsection{Solving the least squares problem}
To start of with our sparse least square problem, we tried LSQR\cite{paige1982lsqr} and LSMR\cite{fong2011lsmr}. However, our problem is extremely ill-conditioned and the residual of both the methods barely converged. Since, SVDs are easier approach to tackle ill-conditioning we move our focus to Krylov based sparse SVD. 

\begin{table}[!h]
\begin{center}
\begin{tabular}{ l | c | c }
  & \textbf{PIELM} & \textbf{Sparse-PIELM (Ours)}\vspace{0.1em}\\
 \hline
 \textbf{Matrix shape} & 5000 $\times$ 1001 & 5000 $\times$ 1001 \vspace{0.1em}\\ 
 \textbf{Matrix rank} & 11 & 561 \vspace{0.1em}\\  
 \textbf{Matrix condition number} & 5.78e+19 & 1.89e+21  \vspace{0.1em} \\
  \textbf{Density} & 1.0 & 0.05  \vspace{0.1em} \\ 
 \textbf{Characteristics} & 1. Dense & 1. Sparse \\
 & 2. Low Rank & 2. Better Controlled Rank\\
 & 3. Extremely Ill-conditioned & 3. Extremely Ill-conditioned \\
 \hline
\end{tabular}
\end{center}
\caption{Characteristics of activation matrices formed in dense (typical) vs our formulation.}
\label{tab:matrices}
\end{table}

\newpage
\subsection{Comparison of Lanczos-Golub-Kahan with and without Full-Orthogonalization for Sparse SVD}

\subsubsection{Case - I: Hard Random Matrix}
\textbf{Goal:} Construct a sparse matrix with slow SVD convergence.

\textbf{Low-rank structure:}
\[
A_{\text{low}} = U S V^\top, \quad
U \in \mathbb{R}^{m \times r},\; V \in \mathbb{R}^{n \times r}
\]
\[
S = \operatorname{diag}(\sigma_i), \quad \sigma_i = i^{-1/2}
\]

\textbf{Sparse noise:}
\[
(A_{\text{sparse}})_{ij} =
\begin{cases}
\mathcal{N}(0,1), & \text{w.p. } \rho \\
0, & \text{otherwise}
\end{cases}
\]

\textbf{Final matrix:}
\[
A = A_{\text{low}} + \varepsilon A_{\text{sparse}}
\]

\textbf{Why difficult?}
\begin{itemize}
    \item Slow decay $\sigma_i \sim i^{-1/2}$ $\Rightarrow$ small spectral gaps
    \item Sparse noise $\Rightarrow$ disrupts Krylov subspace alignment
    \item $\Rightarrow$ slow convergence of Lanczos / Golub--Kahan
\end{itemize}

\begin{figure}[H]
    \centering
    \includegraphics[width=\linewidth]{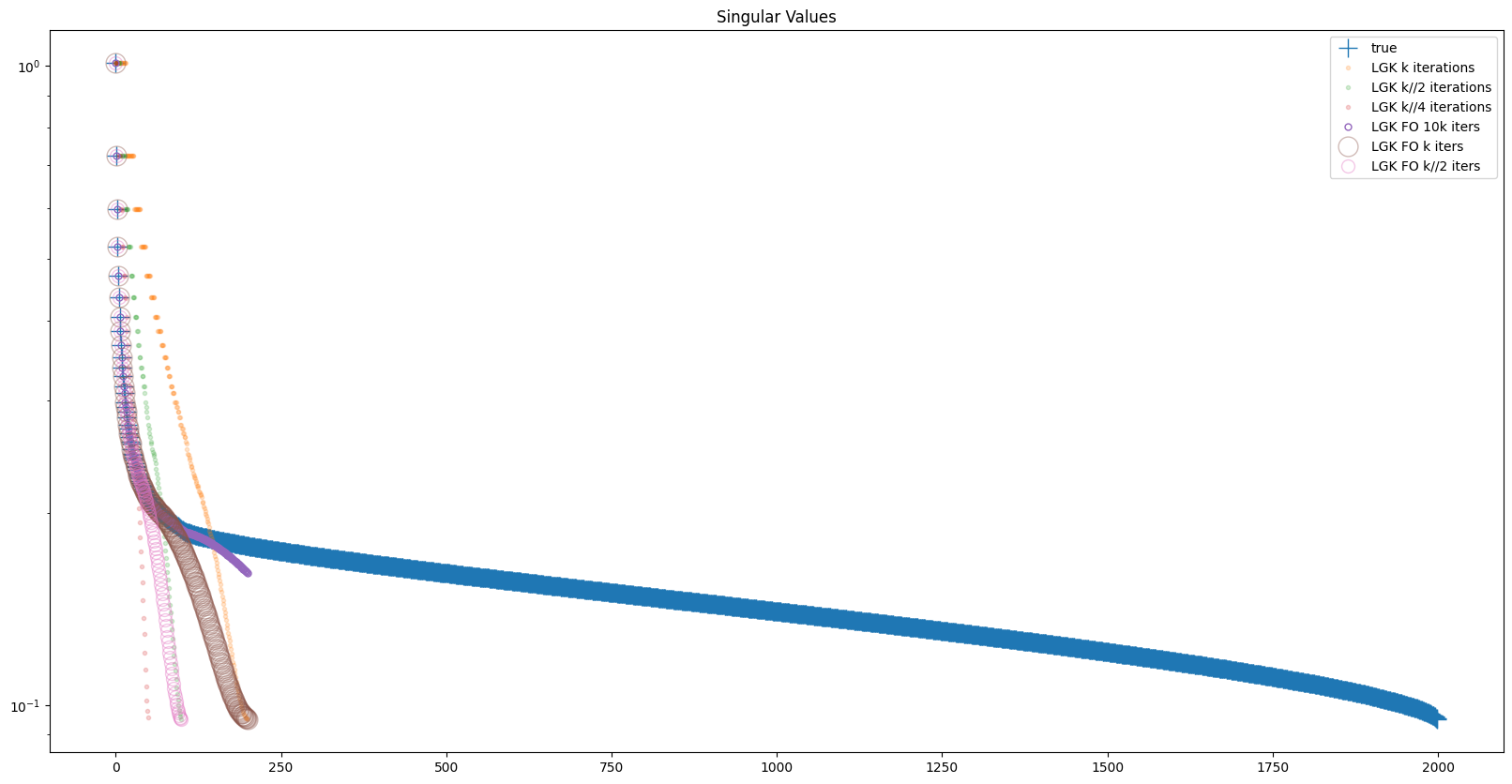}
    \caption{Singular values for Sparse Hard random matrix $10000 \times 2000$ truncated at 200 values}
    \label{fig:hardrandom3}
\end{figure}

\begin{figure}[H]
    \centering
    \includegraphics[width=\linewidth]{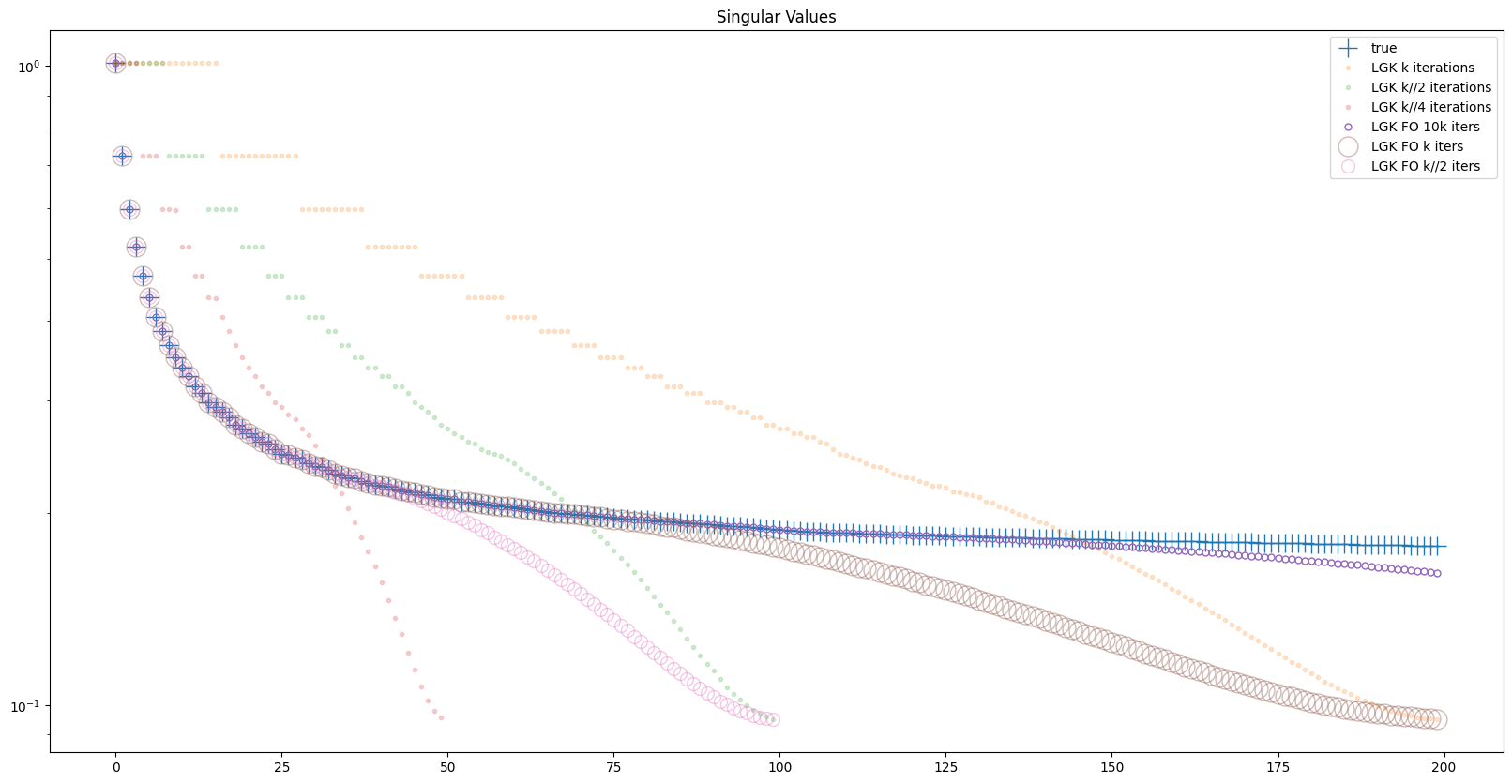}
    \caption{Degenerate singular values for Sparse Hard random matrix $10,000 \times 2,000$ without full-orthogonalization}
    \label{fig:hardrandom2}
\end{figure}

\begin{figure}[H]
    \centering
    \begin{subfigure}{0.48\linewidth}
        \centering
        \includegraphics[width=\linewidth]{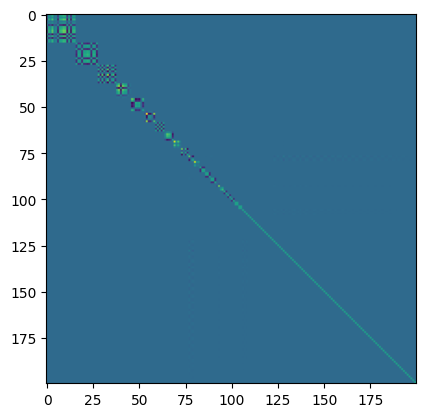}
        \caption{$U^\top U$ with standard bidiagonalization}
    \end{subfigure}
    \hfill
    \begin{subfigure}{0.48\linewidth}
        \centering
        \includegraphics[width=\linewidth]{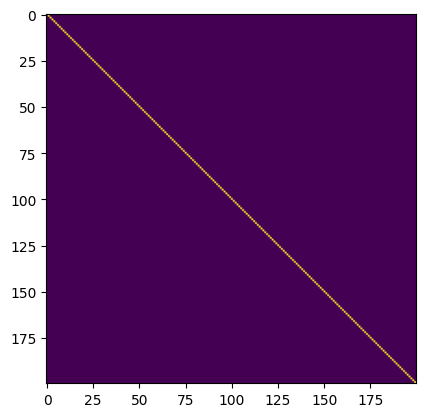}
        \caption{$U^\top U$ with bidiagonalization with full orthogonalization}
    \end{subfigure}
    \caption{Visualization of supposed to be orthonormal matrices.}
    \label{fig:orthogonalityComparison}
\end{figure}

Below in \Cref{fig:hardrandom2,fig:hardrandom3}, we present the results of the bidiagonalization algorithm run on a randomly generated sparse matrix of size $10,000 \times 2,000$. Without orthogonalization, the singular values cover the range but they are inexact. But most specially, it should be observed that without full orthogonalization the produced singular values are degenerate (almost duplicate), which increases with increasing number of Lanczos iterations. \Cref{fig:orthogonalityComparison} shows that without full orthogonalization, the $U$ and $V$ matrices are not purely orthonormal, hence the $U^TU$ is clearly not an identity matrix in the case without full-orthogonalization. The true singular values are obtained by direct method.

\begin{figure}[H]
    \centering
    \includegraphics[width=\linewidth]{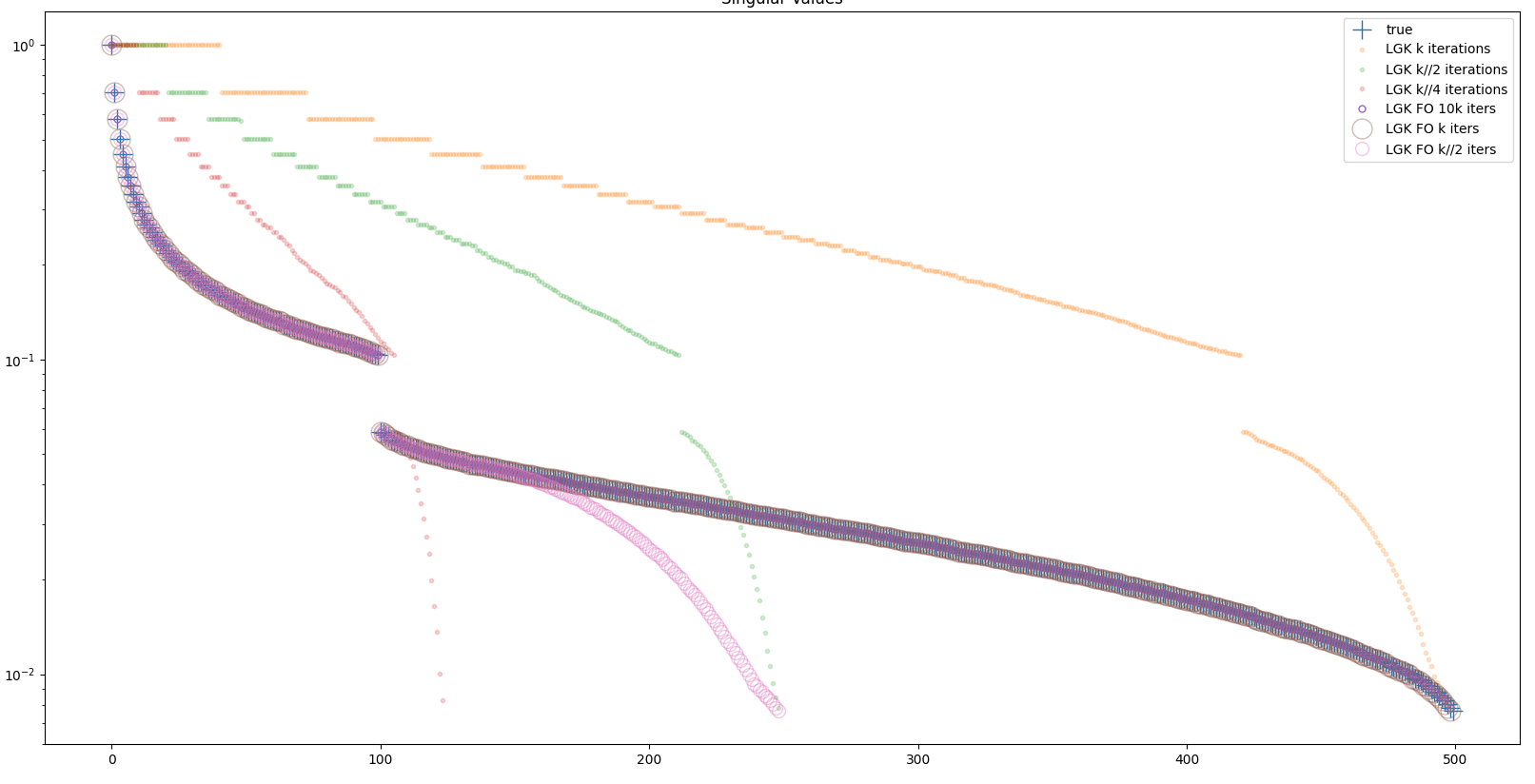}
    \caption{Singular values for Sparse Hard random matrix $2,000 \times 500$ with rank $100$.}
    \label{fig:HardLowRank}
\end{figure}

Similarly, we take a sparse overdetermined matrix of size $20,000 \times 500$ with a rank equal to $100$. The obtained singular values can be seen in \Cref{fig:HardLowRank}. Another important observation is that the range of primary singular values and the singular values of the noise are captured decently but to identify the primary singular values any decent number of iterations is sufficient with full-orthogonalization. However, to exactly capture the non-dominant singular values a very large number of iterations is needed.

\subsubsection{Case - II: Brazilian Interconnect Power System Xingo (CEPEL, Brazil) Extremely Ill-conditioned}

\begin{figure}[H]
    \centering
    \includegraphics[width=\linewidth]{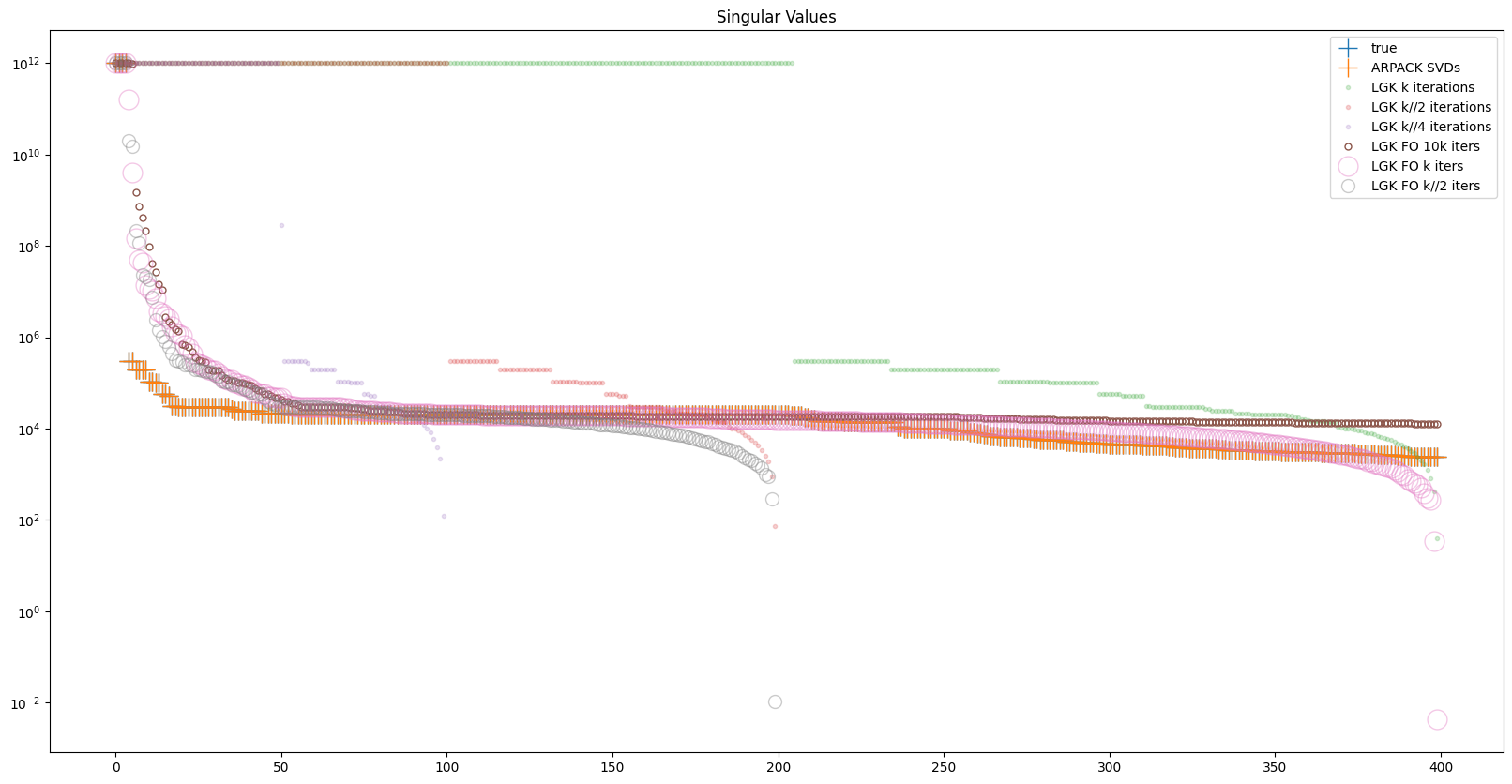}
    \caption{Singular values obtained using sparse SVD using LGK with and without full-orthogonalization.}
    \label{fig:Xingo1}
\end{figure}

\begin{figure}[H]
    \centering
    \includegraphics[width=\linewidth]{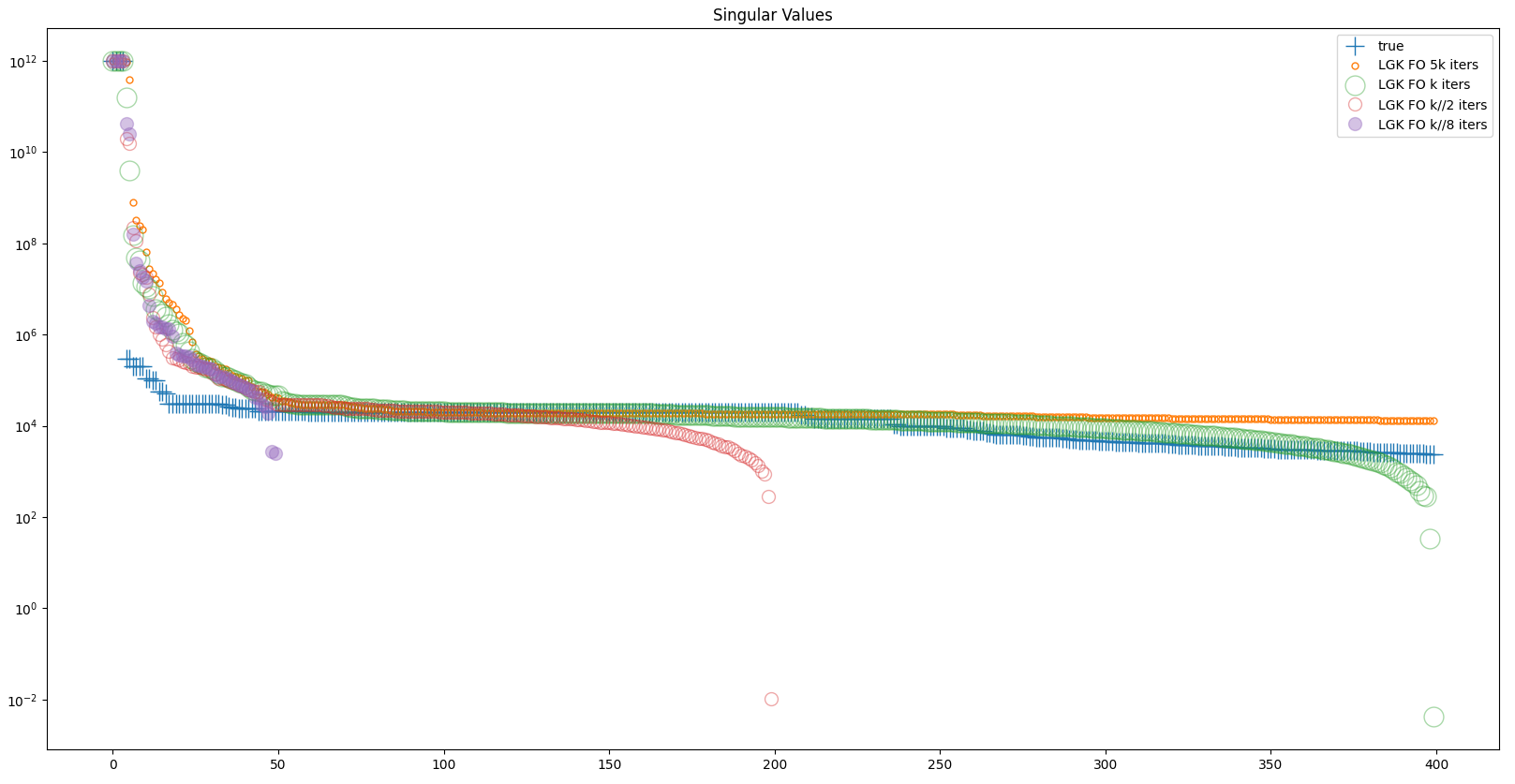}
    \caption{Singular values obtained using sparse SVD using LGK with full-orthogonalization.}
    \label{fig:Xingo2}
\end{figure}

We take a benchmark ill-conditioned sparse matrix from Brazilian Interconnect Power System called Xingo3012. \Cref{tab:xingo3012_info} shows the statistical details of the table. We present the singular values obtained using sparse SVD using Lanczos-Golub-Kahan iterations and full-orthogonalization using Gram-Schmit in \Cref{fig:Xingo1,fig:Xingo2}. If $k$ is the truncation index we need minimum $k$-iterations to get $k$ to obtain all the $k$ singular values. Without full-orthogonalization it is generates a lot of degenerate values even with a small fraction of iterations. While with full-orthogonalization 10$k$ iterations may not be sufficient as we can see in the figures. The range is captured decently well but the intermediate values are very far from actual values. Except the first singular value (maximum), all our solutions are far from actual initial set of singular values. Without full-orthogonalization we rather see an advantage, here, the singular values are approximately in the range of true singular values. With full-orthogonalization we get see see many incorrect singular values out of range, specifically when in between the first and second maximum true singular values.  

\begin{table}[H]
\centering
\begin{tabular}{|c|c|}
\hline
\textbf{Metric} & \textbf{Value} \\
\hline
Dimensions $(m \times n)$ & 20,944 $\times$ 20,944 \\
\hline
Number of Nonzeros (nnz) & 74,386 \\
\hline
Density & 0.0169 \% \\
\hline
Minimum Singular Value &	5.65e-23 \\
\hline
Condition Number &	\textbf{1.77e+34} \\
\hline
Rank (effective / numerical) & 8,750/20,943 \\
\hline
\end{tabular}
\vspace{0.5em}
\caption{Structural and spectral properties of the Xingo3012 sparse matrix.}
\label{tab:xingo3012_info}
\end{table}

\newpage

\subsection{Solving Steady 1D Convection-Diffusion using sparse SVD}

\begin{figure}[H]
    \centering
    \begin{subfigure}{0.48\linewidth}
        \centering
        \includegraphics[width=\linewidth]{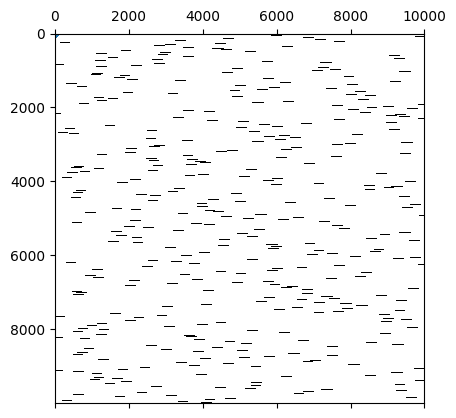}
        \caption{Activation matrix generated with kernel width $1e-5$.\\ Density is $0.029$.}
    \end{subfigure}
    \hfill
    \begin{subfigure}{0.48\linewidth}
        \centering
        \includegraphics[width=\linewidth]{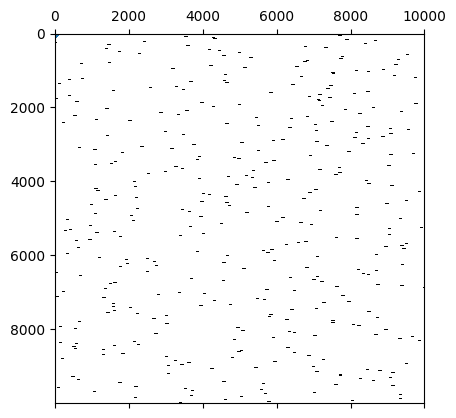}
        \caption{Activation matrix generated with kernel width $1e-6$.\\ Density is $0.009$.}
    \end{subfigure}
    \caption{Spy on activation matrices.}
    \label{fig:densityComparison}
\end{figure}

\begin{figure}[H]
    \centering
    \begin{subfigure}{0.48\linewidth}
        \centering
        \includegraphics[width=\linewidth]{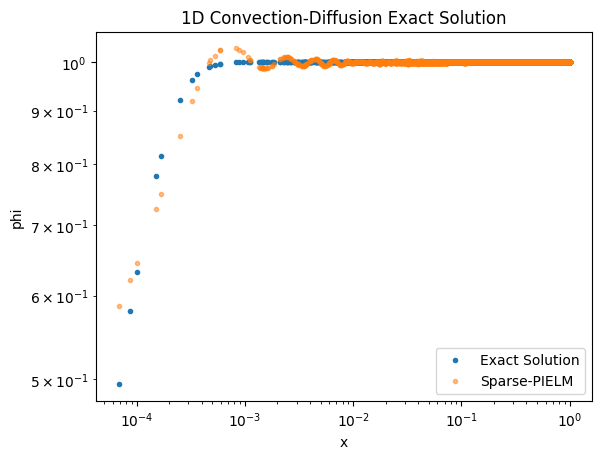}
        \caption{Kernel width $1e-5$.}
    \end{subfigure}
    \hfill
    \begin{subfigure}{0.48\linewidth}
        \centering
        \includegraphics[width=\linewidth]{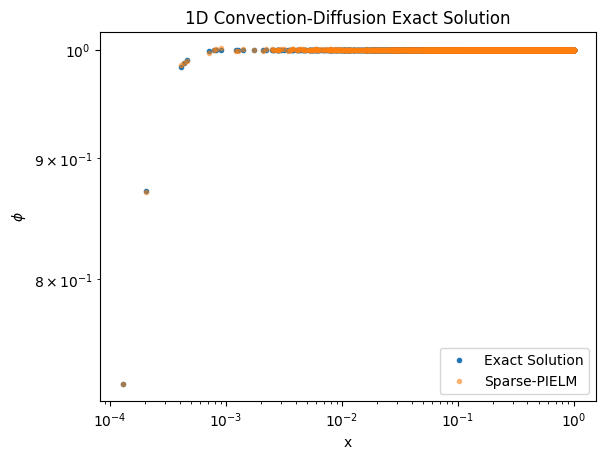}
        \caption{Kernel width $1e-6$,}
    \end{subfigure}
    \caption{Our solution of 1D steady convection-diffusion equation for $Pe = 1e4$.}
    \label{fig:solComparison}
\end{figure}

The PIELM was not improving the results with larger matrices. It could solve up to, at best, $Pe = \{1-5\}e1$. While our sparse formulation scaled well with matrix size, till $Pe = 1e3$ but could not be worked with as it required much larger memory than the resources available in my personal computers. We solved the problem in \Cref{problem} but with $Pe = 1e4$. We took $10,000$ nodes in the network and $10,000$ collocation points which solved $10,000 \times 10,000$ sparse matrix in least square setup using Lanczos-Golub-Kahan Bidiagonalization with full-orthogonalization and restarting truncated at $5,000$ iterations. The activation matrix is generated after encoding using Gaussian kernel with width $1e-5$, the density of the matrix is $0.029$ does not solve our problem. The formed matrix is low rank and hence less expressive. Hence, we are motivated the reduce the kernel width by a factor of $10$. The activation matrix is generated after encoding using Gaussian kernel with width $1e-6$, the density of the matrix is $0.00932$, which is less than $< 1 \%$ and it has a much higher rank. The wall-clock runtime was around 1 hour and 57 minutes. See \Cref{fig:densityComparison} to spy on the density of the two matrices obtained by randomly ordered $x$ values, order $x$ values gives a structured matrix as seen in \Cref{fig:matrices}. The boundary layer solution for the two kernel width setting in log-log scale for the case can be seen in \Cref{fig:solComparison}. By exploiting sparsity we are able to scale and improve the results for a tough problem. 

\textbf{Limitations:} If we refer to \Cref{equ:sparsePseudoINV}, the orthonormal matrices $U$ and $V$ are dense, hence the equation forms a large dense matrix for the closed form solution, which shall require large memory in high/full rank setup. We had to increase the rank to improve the expressivity while high-rank is able to solve better but it gives rise to a dense matrix formation. So, they are competing interests.


\section{Conclusion}
In this work, we presented a sparse and linear formulation for solving a tough differential equation with sharp gradient through least-squares-based nonlinear function approximation within the PIELM framework. While this formulation enables efficient training by avoiding iterative gradient-based optimization, it introduces critical numerical challenges arising from low-rank structure and ill-conditioning of the resulting system matrices. We highlighted the limitations of naïve singular value decomposition approaches based on forming $A^\top A$, which can significantly amplify conditioning issues and degrade numerical stability. To address these challenges, we explored Krylov subspace methods, emphasizing the role of Lanczos iterations for symmetric eigenvalue problems and the more general Golub-Kahan bidiagonalization for computing sparse SVD in large-scale settings. The connection $T_k = B_k^\top B_k$ provides a useful bridge between these formulations, enabling efficient approximation of singular values without explicitly forming dense matrices. Finally, we discussed practical concerns such as loss of orthogonality in iterative processes and the necessity of re-orthogonalization strategies to maintain numerical robustness. Finally, we show solve the problem with a $10\times$ larger convection dominated case.

\section*{Acknowledgement}

We are especially thankful to Prof. Chen Greif for valuable guidance, encouragement, and continuous support. Finally, I am grateful to all those who directly or indirectly contributed to this work.

\bibliography{references}
\bibliographystyle{unsrt}

\newpage
\appendix

\end{document}